# ON TWO OF ERDÖS'S OPEN PROBLEMS

Florentin Smarandache
University of New Mexico
200 College Road
Gallup, NM 87301, USA
E-mail: smarand@unm.edu

**Abstract.**
This short note presents some remarks and conjectures on two open problems proposed by P. Erdös.

**First Problem.**

In one of his books ("Analysis…") Mr. Paul Erdös proposed the following problem:
"The integer $n$ is called a barrier for an arithmetic function $f$ if $m + f(m) \leq n$ for all $m < n$.
Question: Are there infinitely many barriers for $\varepsilon v(n)$, for some $\varepsilon > 0$? Here $v(n)$ denotes the number of distinct prime factors of $n$."

We found some results regarding this question, which results make us to conjecture that there is a finite number of barriers, for all $\varepsilon > 0$.

Let $R(n)$ be the relation: $m + \varepsilon v(m) \leq n, \ \forall m < n$.

**Lemma 1.1.** If $\varepsilon > 1$ there are two barriers only: $n = 1$ and $n = 2$ (which we call trivial barriers).

*Proof.* It is clear for $n = 1$ and $n = 2$ because $v(0) = v(1) = 0$.
Let's consider $n \geq 3$. Then, if $m = n - 1$ we have $m + \varepsilon v(m) \geq n - 1 + \varepsilon > n$, contradiction.

**Lemma 1.2.** There is an infinity of numbers which cannot be barriers for $\varepsilon v(n)$, $\forall \varepsilon > 0$.

*Proof.* Let's consider $s, k \in \mathbb{N}^*$ such that $s \cdot \varepsilon > k$. We write $n$ in the form $n = p_{i_1}^{\alpha_{i_1}} \cdots p_{i_s}^{\alpha_{i_s}} + k$, where for all $j$, $\alpha_{i_j} \in \mathbb{N}^*$ and $p_{i_j}$ are positive distinct primes.

Taking $m = n - k$ we have $m + \varepsilon v(m) = n - k + \varepsilon \cdot s > n$.

But there exists an infinity of $n$'s because the parameters $\alpha_{i_1}, ..., \alpha_{i_s}$ are arbitrary in $\mathbb{N}^*$ and $p_{i_1}, ..., p_{i_s}$ are arbitrary positive distinct primes, also there is an infinity of couples $(s, k)$ for an $\varepsilon > 0$, fixed, with the property $s \cdot \varepsilon > k$.

**Lemma 1.3.** For all $\varepsilon \in (0,1]$ there are nontrivial barriers for $\varepsilon v(n)$.

*Proof.* Let $t$ be the greatest natural number such that $t\varepsilon \leq 1$ (always there is such $t$).

Let $n$ be from $[3,..., p_1 \cdots p_t p_{t+1})$, where $\{p_i\}$ is the sequence of the positive primes. Then $1 \leq v(n) \leq t$.

All $n \in [1,..., p_1 \cdots p_t p_{t+1}]$ is a barrier, because: $\forall\ 1 \leq k \leq n-1$, if $m = n - k$ we have $m + \varepsilon v(m) \leq n - k + \varepsilon \cdot t \leq n$.

Hence, there are at list $p_1 \cdots p_t p_{t+1}$ barriers.

**Corollary.** If $\varepsilon \to 0$ then $n$ (the number of barriers) $\to \infty$.

**Lemma 1.4.** Let's consider $n \in [1,..., p_1 \cdots p_r p_{r+1}]$ and $\varepsilon \in (0,1]$. Then: $n$ is a barrier if and only if $R(n)$ is verified for $m \in \{n-1, n-2, ..., n-r+1\}$.

*Proof.* It is sufficient to prove that $R(n)$ is always verified for $m \leq n - r$.

Let's consider $m = n - r - u$, $u \geq 0$. Then $m + \varepsilon v(m) \leq n - r - u + \varepsilon \cdot r \leq n$.

**Conjecture.**

We note $I_r \in [p_1 \cdots p_r, ..., p_1 \cdots p_r p_{r+1})$. Of course $\bigcup_{r \geq 1} I_r = \mathbb{N} \setminus \{0,1\}$, and $I_{r_1} \cap I_{r_2} = \Phi$ for $r_1 \neq r_2$.

Let $\mathcal{N}_r(1+t)$ be the number of all numbers $n$ from $I_r$ such that $1 \leq v(n) \leq t$.

We conjecture that there is a finite number of barriers for $\varepsilon v(n)$, $\forall \varepsilon > 0$; because

$$\lim_{r \to \infty} \frac{\mathcal{N}_r(1+t)}{p_1 \cdots p_{r+1} - p_1 \cdots p_r} = 0$$

and the probability (of finding of $r-1$ consecutive values for $m$, which verify the relation $R(n)$) approaches zero.

**Second Problem.**

Paul Erdös has proposed another problem:
(1) "Is it true that $\lim_{n \to \infty} \max_{m < n}(m + d(m)) - n = \infty$?, where $d(m)$ represents the number of all positive divisors of $m$."

We clearly have :

**Lemma 2.1.** $(\forall) n \in \mathbb{N} \setminus \{0,1,2\}$, $(\exists)! s \in \mathbb{N}^*$, $(\exists)! \alpha_1,..., \alpha_s \in \mathbb{N}$, $\alpha_s \neq 0$, such that $n = p_1^{\alpha_1} \cdots p_s^{\alpha_s} + 1$, where $p_1, p_2, ...$ constitute the increasing sequence of all positive primes.

**Lemma 2.2.** Let $s \in \mathbb{N}^*$. We define the subsequence $n_s(i) = p_1^{\alpha_1} \cdots p_s^{\alpha_s} + 1$, where $\alpha_1, ..., \alpha_s$ are arbitrary elements of $\mathbb{N}$, such that $\alpha_s \neq 0$ and $\alpha_1 + ... + \alpha_s \to \infty$ and we order it such that $n_s(1) < n_s(2) < ...$ (increasing sequence).

We find an infinite number of subsequences $\{n_s(i)\}$, when $s$ traverses $\mathbb{N}^*$, with the properties:

a) $\lim_{i \to \infty} n_s(i) = \infty$ for all $s \in \mathbb{N}^*$.

b) $\{n_{s_1}(i), i \in \mathbb{N}^*\} \cap \{n_{s_2}(j), j \in \mathbb{N}^*\} = \Phi$, for $s_1 \neq s_2$ (distinct subsequences).

c) $\mathbb{N} \setminus \{0,1,2\} = \bigcup_{s \in \mathbb{N}^*} \{n_s(i), i \in \mathbb{N}^*\}$

Then:

**Lemma 2.3.** If in (1) we calculate the limit for each subsequence $\{n_s(i)\}$ we obtain:

$$\lim_{n \to \infty} \left( \max_{m < p_1^{\alpha_1} \cdots p_s^{\alpha_s}} (m + d(m)) - p_1^{\alpha_1} \cdots p_s^{\alpha_s} - 1 \right) \geq \lim_{n \to \infty} \left( p_1^{\alpha_1} \cdots p_s^{\alpha_s} + (\alpha_1 + 1)\ldots(\alpha_s + 1) - p_1^{\alpha_1} \cdots p_s^{\alpha_s} - 1 \right) =$$

$$= \lim_{n \to \infty} \left( (\alpha_1 + 1)\ldots(\alpha_s + 1) - 1 \right) > \lim_{n \to \infty} (\alpha_1 + \ldots + \alpha_s) = \infty$$

From these lemmas it results the following:

**Theorem:** We have $\overline{\lim_{n \to \infty}} \max_{m < n} (m + d(m)) - n = \infty$.

**REFERENCES**

[1] P. Erdös - Some Unconventional Problems in Number Theory - Mathematics Magazine, Vol. 57, No.2, March 1979.

[2] P. Erdös - Letter to the Author - 1986: 01: 12.

[Published in "Gamma", XXV, Year VIII, No. 3, June 1986, p. 5.]